\documentclass[10.5pt]{article}
\usepackage{amsmath}
\usepackage{amsfonts}
\usepackage{amsthm}
\usepackage{amscd}
\usepackage{amssymb}
\usepackage{color}
\usepackage{enumerate}

\newtheorem{theorem}{Theorem}[section]

\newtheorem{corollary}[theorem]{Corollary}

\theoremstyle{definition}

\newtheorem{example}[theorem]{Example}
\theoremstyle{remark}
\newtheorem{remark}[theorem]{Remark}

\newcommand{\1}{\mathbf{1}}

\def\1{\mathbf{1}}

\def\bbZ{\mathbb{Z}}

\def\f2{\mathbb{F}_2}

\newcommand{\bbN}{\mathbb{N}}
\newcommand{\al}{\alpha}

\newcommand{\de}{\delta}
\newcommand{\hde}{{\tilde\de}}
\newcommand{\e}{\varepsilon}

\newcommand{\la}{\lambda}

\newcommand{\s}{\sigma}

\newcommand{\bbR}{\mathbb{R}}

\newcommand{\Span}{\operatorname{span}}
\newcommand{\codim}{\operatorname{codim}}

\newcommand{\lb}{\label}

\newcommand{\wtw}{if and only if}

\newcommand{\Buo}{Without loss of generality }

\newcommand{\DEF}{\buildrel {\mbox{\tiny def}}\over =}

\newcommand\remove[1]{}

\usepackage{mathtools,graphicx}
\newcommand\vertarrowbox[2]{%
    \begin{array}[t]{@{}c@{}} #1 \\
    \rotatebox{90}{$\xrightarrow{\hphantom{}}$} \\[-1.5ex]
    \mathclap{\scriptstyle\text{#2}}%
    \end{array}}

\numberwithin{equation}{section}

\begin{document}

\title{On  the relations between Auerbach or almost Auberbach Markushevich systems and Schauder bases}

\author{Beata Randrianantoanina, Micha\l \ Wojciechowski\thanks{M.W. supported by the National Science Centre, Poland, CEUS - UNISONO programme, 
project no. 2020/02/Y/ST1/00072},  and Pavel Zatitskii}

\date{~}

\maketitle

\begin{abstract}
We establish that the summability of the series $\sum\varepsilon_n$ is the necessary and sufficient criterion ensuring that every
$(1+\varepsilon_n)$-bounded Markushevich basis in a separable Hilbert space is a Riesz basis. Further we show that if $n\varepsilon_n\to \infty$, then in $\ell_2$ there exists a $(1+\varepsilon_n)$-bounded 
 Markushevich basis that under any permutation is  non-equivalent to a Schauder basis. We extend this result to any separable Banach space. Finally we provide examples of Auerbach bases in  
 1-symmetric separable Banach  spaces whose no permutations are equivalent to any Schauder basis or (depending on the space) any unconditional Schauder basis.
\end{abstract}

{\small \noindent{\bf Keywords:} Schauder basis; Markushevich basis; Auerbach basis;

}

\section{Introduction}

One of the main questions at the early stage of the development of the theory  of infinite dimensional linear spaces was the question of  existence of reasonable coordinates which would allow the use of the Cartesian method analogously as in finite dimensional vector spaces.
Parseval's theorem solves this problem completely in the most prominent case of Hilbert spaces. However  the case of general Banach spaces  turned out to be much more complicated.  
It is clear that individual coordinates of each point should arise as evaluations of  appropriately chosen universal set of
 linear functionals. However the problem of identifying  ``good''  systems of linear functionals   remained open for a long time. 

An early natural candidate for such a system was introduced by Schauder in 1927 \cite{S27},   who defined a basis of a Banach space, now called {\it a Schauder basis}, to be a system of  biorthogonal sequences of vectors and linear functionals, such that every element of the space is uniquely represented by a limit of   uniformly bounded partial sums, see e.g. \cite[Definition 1.a.1]{LT1}. 
It has been known since 1930s that  Schauder bases exist in important known examples of classical Banach spaces, such as Lebesgue spaces $L_p$, spaces of continuous functions $C(K)$, 
Hardy spaces $H^p$, Sobolev spaces $W_{p.k}$. The problem whether   every separable Banach space has a Schauder basis has appeared already in Banach's book \cite{Banach} in 1932, but despite a great deal of research by many mathematicians, it remained open for almost forty years. It was not until after Grothendieck  undertook the study of variants of the notion of an approximation property in the fifties \cite{G55}, that,  in the beginning of 1970s, Enflo \cite{E73} constructed first  examples of separable Banach spaces without a Schauder basis. 

Since the problem of the existence of a Schauder basis turned out to be very complicated, it was natural to consider  possible weaker notions, in particular, what happens if  the requirement 
of uniform boundedness of partial sums was relaxed. 
Answering this weaker question,  in 1943, Markuschevich \cite{M43}  proved that  every separable Banach space $X$ contains a biorthogonal  sequence of vectors and   functionals
 $\{x_n,x_n^*\}_{n\in \bbN}$ in $X\times X^*$, such that   $\{x_n\}_{n\in \bbN}$ is {\it fundamental}, i.e. the set of all linear combinations of $x_i$'s is  norm dense in $X$, and
  $\{x_n^*\}_{n\in \bbN}$ is {\it total}, i.e. $x=0$ if  $x_n^*(x)=0$ for all $n\in \bbN$. Such a biorthogonal system is now  called a {\it Markushevich basis}, or an {\it M-basis} for short.

It was not until over thirty years later, that, in 1975, Ovsepian and Pe\l czy\'nski \cite{OP75},  proved that in any separable Banach space there exists a Markushevich basis for which both the  sets of  vectors and of  functionals are bounded.
While  the above result has a  rather isomorphic flavor,
in 1976, Pe\l czy\'nski \cite{P76} and  Plichko \cite{Pl76},  independently, using Dvoretzky theorem, proved that  the bound may be requested to be  arbitrarily close to 1 (cf. also  \cite{Pl77,Pl80}). After another 20 years - in  1999 - Vershinin \cite{V99} further refined Pe\l czy\'nski's argument and proved  that in every separable Banach space $X$, for every nonnegative sequence 
$(\e_n)_n$ with $\sum\varepsilon_n^2=\infty$, there exists an M-basis 
  $\{x_n,x_n^*\}_{n\in \bbN} \subseteq X\times X^\ast $ such that $\|x_n\|\cdot\|x_n^\ast\|\leq 1+\varepsilon_n$ (such M-bases are called {\it $(1+\e_n)$-bounded}).

The ideal situation  would be if one could take $\varepsilon_n=0$ for every $n$. In 1930 in his PhD thesis \cite{A1930} (see \cite[Proposition 1.c.3]{LT1})  Auerbach  proved  that this is possible   in every Banach space of finite dimension.  
In 2017, answering a question of Pe\l czy\'nski, Weber and Wojciechowski \cite{WW2017} proved that in every $n$ dimensional Banach space this is possible to do it in at least
$\frac12n(n-1)+1$ substantially different ways.

An M-basis such that  $\|x_n\|=\|x_n^\ast\|= 1$, for all $n$, is called 
an {\it Auerbach basis}. A $(1+\e_n)$-bounded M-basis is called an {\it almost Auerbach basis}.

 In  infinite dimensional separable Banach spaces, apart from some very special constructions, our knowledge about  Auerbach bases is highly unsatisfactory. In particular, Banach's question from 1932 \cite[Remarques. Chapitre VII. \S1, p. 237-238]{Banach} whether  every separable Banach space has an Auerbach basis still remains open.  
 
 On the one hand, 
 Plichko  \cite[Corollary 1]{Pl77}
 proved  that every separable  Banach space $X$ for every $\e>0$  has a  $(1+\e)$-equivalent norm
  $|\!|\!|\cdot|\!|\!|$ such that $(X,|\!|\!|\cdot|\!|\!|)$ has an Auerbach basis (cf.  \cite[Proposition 1.31]{HMVZ}).
 On the other hand, Vershinin's  result \cite{V99}  mentioned above, states that every separable Banach space $X$ has an almost Auerbach basis, provided that the sequence $(\e_n)_n$ is not 
 square summable -  till now, this is  the 
 most far-reaching result towards establishing the  existence of an Auerbach basis in a given Banach space.

Since the  notions of  Markushevich and Auerbach bases are formally weaker than that of a Schauder basis, it is natural to ask  what conditions on  a Markushevich  basis  would guarantee that it is also a Schauder basis. 

Clearly, in a Hilbert space any Auerbach basis   is an orthonormal basis.
However,  in general spaces, as examples of trigonometric systems in $L^1(T)$ and $C(T)$  show, even being an Auerbach basis does not ensure being  a Schauder basis  (see e.g. \cite[Section 1.f]{LT1}). Moreover, Johnson proved that every separable Banach space contains a (usually unbounded) 
M-basis that is not a Schauder basis under any permutation, see \cite[Remark 1.33]{HMVZ}. We improve this result in Corollary~\ref{(5)}.

In the present  paper we study  conditions on almost Auerbach M-bases to be Schauder bases and we provide several constructions of  almost Auerbach and Auerbach bases in separable Banach spaces that, under any permutation, are not equivalent to a Schauder basis.

We start by proving a certain
stability property of Auerbach bases in Hilbert spaces. Specifically, we identify a necessary and sufficient condition for a sequence 
$(\varepsilon_n)_n$ so that in a Hilbert space every   $(1+\e_n)$-bounded M-basis is equivalent to an orthonormal basis.

\begin{theorem}\lb{(1)}
For any     sequence of nonnegative numbers $(\e_n)_{n=1}^\infty$ such that $ \sum_{n=1}^\infty \e_n<\infty$, every normalized $(1+\e_n)$-bounded M-basis of $\ell_2$  is  equivalent to an orthonormal   basis of $\ell_2$. 
\end{theorem}

Recall that two  minimal  systems 
$\{x_n\}_n$ and $\{y_n\}_n$    
are called {\it equivalent} \wtw\ the linear map $A$  defined by $Ax_n=y_n$, for all $n$, is a bounded 
isomorphism of 
$\overline{\Span}\{x_n\}_n$ onto $\overline{\Span}\{y_n\}_n$, that is,  $\|A\|\cdot\|A^{-1}\|<\infty$ (cf. \cite[Section 1.f]{LT1}). In this case, the number $\|A\|\cdot\|A^{-1}\|$ is called the {\it distance} between these
two  systems.

\begin{theorem}\lb{(2)}
 For any sequence of nonnegative numbers $(\e_n)_{n=1}^\infty$ such that $ \sum_{n=1}^\infty \e_n=\infty$, there exists an $(1+\e_n)$-bounded M-basis of 
$\ell_2$ which is not equivalent to an orthonormal basis of $\ell_2$. 
\end{theorem}

 The necessary condition could be derived from the known constructions of conditional bases in a Hilbert space, but  we also provide a necessary condition on 
$(\varepsilon_n)$ for the existence of a $(1+\varepsilon_n)$-M-basis that is not  equivalent to a Schauder basis under any permutation - first in a Hilbert space and then in  an arbitrary separable Banach space.

\begin{theorem}\lb{(4)}
 For every sequence of positive numbers $(\e_n)_{n=1}^\infty$  such that 
\begin{equation}\lb{extra}
\lim_{n\to\infty}n\e_n=\infty,
\end{equation}
 there exists an $(1+\e_n)$-bounded M-basis of $\ell_2$, such that none of its permutations is    equivalent to  a Schauder basis of $\ell_2$. 
\end{theorem}

Note that \eqref{extra} clearly implies that $ \sum_{n=1}^\infty \e_n=\infty$.

As a consequence of Theorem~\ref{(4)}, the result of Vershynin \cite{V99}, and the Dvoretzky Theorem,
we obtain:

\begin{corollary}\lb{(5)}
In  any separable  Banach space   $X$, for any     sequence of positive numbers  
$(\e_n)_{n\in\bbN}$ that satisfies  \eqref{extra} and such that 
 $\sum \e_n^2=\infty$, there exists a
  $(1+\e_n)$-bounded M-basis   such that none of its permutations is    equivalent to  a Schauder basis of $X$. 
\end{corollary}

The next result  follows from Corollary~\ref{(5)} using the same reasoning as Plichko in the proof of \cite[Corollary 1]{Pl77}.

\begin{corollary}\lb{(3)}
For every separable  Banach space   $(X,\|\cdot\|)$ and  every $\e>0$,  there exists an equivalent norm $|\!|\!|\cdot|\!|\!|$ on $X$ such that
  $\|x\|\le |\!|\!|x|\!|\!|\le (1+\e)\|x\|$ and $(X,|\!|\!|\cdot|\!|\!|)$ has an Auerbach basis
  such that none of its permutations is    equivalent to  a Schauder basis of $X$. 
\end{corollary}

\begin{remark}
It also follows from the proof of Corollary~\ref{(5)} that, for any sequence $(\e_n)_{n\in\bbN}$ that is not square summable and satisfies  \eqref{extra},  one could additionally request that the equivalent norm $|\!|\!|\cdot|\!|\!|$ on $X$  in Corollary~\ref{(3)}  is such that for some nested sequence of finite-codimensional subspaces $X_k$ of $X$ (with $\codim X_k$ growing fast enough) we have 
 $\|x\|\le |\!|\!|x|\!|\!|\le (1+\e_k)\|x\|$ for all $x\in X_k$ and all $k$.
\end{remark}

We finish the paper with some examples. In Example~\ref{exl1}, using recent deep results from  combinatorial number theory \cite{gs,Sa}, we construct in the space 
$\ell_1$ an Auerbach basis that is not equivalent to a Schauder basis under any permutation.
In Example~\ref{exlp} we show that
for any  $1<p<\infty$, $p\neq 2$, the space $\ell^p$ contains an Auerbach basis that is not an unconditional Schauder basis, and in Example~\ref{exsym} we extend this construction to a wider class of 1-symmetric sequence spaces not isomorphic to a Hilbert space. 

We do not know whether Auerbach bases that are not equivalent to Schauder bases under any permutation  exist  in all  1-symmetric sequence spaces not isomorphic to $\ell_2$, 
in particular, in $\ell^p$, $1<p<\infty$, $p\neq 2$. As we describe in Example~\ref{expel}, this question is related to a 
conjecture of   Pe\l czy\'nski from 2006 \cite[Problem 5.1]{Pe2006}.

In view of our results it is natural to ask whether or not the Hilbert space is the only separable Banach space for which every Auerbach basis  is  equivalent to an unconditional Schauder basis or even  merely to a  Schauder basis.

Theorem~\ref{(1)} suggests a stronger version of the above question, namely whether there exists a separable Banach space other than the Hilbert space such that all sufficiently tightly bounded almost Auerbach bases  have to be  equivalent to an unconditional  basis, or merely  to a  Schauder basis.

Note that Corollary~\ref{(3)} demonstrates that all these questions are strictly  isometric in nature.
 
In the paper we work mainly in the Hilbert space setting  and the results for arbitrary Banach spaces are derived from similar ones for the Hilbert space in a similar way as in
Pełczy\'nski or Vershinin proofs - using Dvoretzky theorem.
For the Hilbert space case, in Section~\ref{prelim} we reduce the construction of a suitable 
M-basis to  spectral analysis of selfadjoint operators, which allows us in subsequent sections to reduce  the proofs of our results to an appropriate  choice of aforementioned spectral parameters. 
Note that our constructions of  ``pathological'' almost Auerbach bases in  Hilbert spaces use finite dimensional arguments which are then glued using  an  orthogonal finite dimensional decomposition of $\ell_2$. We note that, as a consequence, all almost Auerbach bases that we construct are {\it strong}, i.e. for all $x\in X$, 
$x\in\overline{\Span}\{\langle x,x^*_n\rangle x_n\}_{n\in\bbN}$, cf. \cite[Section~1.5]{HMVZ}.

\section*{Acknowledgements}
The authors thank Karol Leśnik and Maciej Rzeszut for discussions during the preparation of this paper, and Anton Tselishchev and the anonymous referee for their questions that led to an improvment of the presentation of the proof of Theorem~\ref{(4)}.

\section{Preliminaries}\lb{prelim}
 
 We use standard definitions and notation, for undefined terms  we refer the reader to \cite{LT1} or \cite{HMVZ}. 

Let $X$ be a Banach space.
A set  of vectors  $\{x_i\}_{i\in I}\subset X$ is called {\it minimal} if there exists a set  
$\{x_i^*\}_{i\in I}\subset X^*$ such that the system  $\{x_i,x_i^*\}_{i\in I}$ is biorthogonal.

A biorthogonal system $\{x_i,x_i^*\}_{i\in I}$ in $X\times X^*$ is called {\it semi-normalized} if there exist constants $c_1,c_2>0$ such that $c_1\le \|x_i\|\le c_2$ and $c_1\le \|x_i^*\|\le c_2$ for each $i\in I$.

  If $\e>0$, we say that an M-basis 
 $\{x_n,x_n^*\}_{n\in \bbN}$ is  {\it $(1+\e)$-bounded} if $\|x_n\|\cdot\|x_n^*\|\le 1+\e$ for all $n\in\bbN$. 
 Given   a sequence of nonnegative numbers $(\e_n)_{n\in\bbN}$,
 we  say that an M-basis $\{x_n,x_n^*\}_{n\in \bbN}$ is {\it $(1+\e_n)$-bounded}, if $\|x_n\|\cdot\|x_n^*\|\le 1+\e_n$ for all $n\in\bbN$. An M-basis
  $\{x_n,x_n^*\}_{n\in \bbN}$ is called an {\it Auerbach basis} if $\|x_n\|=\|x_n^*\|=1$ for all $n\in\bbN$. 
  An  $(1+\e_n)$-bounded M-basis is   called {\it an  almost Auerbach basis}.

 The goal of this section is to describe a general form of bi-orthogonal systems in finite dimensional Hilbert spaces that we will use throughout this paper, and to give the necessary and sufficient quantitative conditions for an M-basis in this form to be almost Auerbach.
 
 Fix $n\in \bbN$,  an 
orthonormal basis 
$\{e_1,\dots,e_n\}$ in $\ell_2^n$, and its dual basis $\{e_1^*,\dots,e_n^*\}$.

Let $(\e_i)_{i=1}^n$ be a sequence of nonnegative numbers
and $\{x_i,x_i^*\}_{i=1}^n$ be a 
$(1+\e_i)$-bounded M-basis in $\ell_2^n$. \Buo we may and do assume throughout this paper that in all considered sequences we have   
$\e_i<1$ for all $i$. 

Let $A$ be the $n\times n$ matrix such that $Ae_i = x_i$. The  columns of $A$
are the vectors $x_1,\dots,x_n$ 
represented in the basis $\{e_i\}$. Let $B=A^TA$.
Then 
\begin{equation}\lb{bii}
b_{ii}=(A^TA)_{ii}=\langle x_i,x_i \rangle=\|x_i\|^2.
\end{equation}

Since the system $\{x_i,x_i^*\}$ is biorthogonal, the vectors $\{x_1^*,\dots,x_n^*\}$ represented in the basis $\{e_i^*\}$ are the rows of the matrix $(A^{-1})^T$. Thus $\|x_i^*\|^2=(A^{-1}(A^{-1})^T)_{ii}$, that is 
\begin{equation}\lb{inverse}
(B^{-1})_{ii}=\|x_i^*\|^2.
\end{equation}

Since the matrix $B$ is positive definite, 
 there exists a unitary matrix $U=[u_{ij}]$ and a diagonal matrix $D=[d_i]=
  \left[ {\begin{array}{ccc}
    d_1 & 0 & 0\\
   0 & \ddots & 0\\
  0   &   0 & d_n
  \end{array} } \right]$,
  with $d_1\ge d_2\ge\dots\ge d_n>0$, 
  such that $B=UDU^T$. Note that by the singular value decomposition there exists a 
  unitary matrix $V$, such that $A=V\sqrt{D}U^T$.
  
   We have
 \begin{align}
 b_{ii}&=\sum_{j=1}^n u_{ij}^2 d_j,\\
   (B^{-1})_{ii}&=\sum_{j=1}^n u_{ij}^2 \frac{1}{d_j}.
   \end{align}
   
   Thus, by \eqref{bii} and \eqref{inverse}, given  a  unitary matrix $ U=[u_{ij}]\in {\mathcal{O}}(n)$, a positive diagonal matrix $D=[d_i]$,
  with $d_1\ge d_2\ge\dots\ge d_n>0$, and  an arbitrary $n\times n$ unitary matrix $V$,  the column vectors of the  the
  matrix $A=V\sqrt{D}U^T$ form a  $(1+\e_i)$-bounded M-basis in $\ell_2^n$   \wtw \ 
\begin{equation}\lb{M-basis-matrix}
\forall 1\le i \le n\ \ \  \ \ \ \left(\sum_{j=1}^n u_{ij}^2 d_j\right)\left(\sum_{j=1}^n u_{ij}^2 \frac{1}{d_j}\right) \le (1+\e_i)^2.
\end{equation}

Note that, since   $A=V\sqrt{D}U^T$, for some $U,V \in {\mathcal{O}}(n)$, and since operators defined by  unitary matrices    are isometries of $\ell_2^n$,  we immediately obtain that the  distance between this M-basis and the orthonormal basis is equal to
\begin{equation}\lb{normofA}
  \|A\|\cdot\|A^{-1}\|= \|\sqrt{D}\|\cdot\|\sqrt{D^{-1}}\|= \frac{\sqrt{d_1}}{\sqrt{d_n}}.
  \end{equation}

\section{Proof of Theorem~\ref{(1)}.}
Let $(\e_i)_{i\in\bbN}$ be a nonnegative sequence  such that  $\e_i<1$, for all $i\in\bbN$,  and
$\sum_{i=1}^\infty \e_i < \infty.$
Let $\{x_i,x_i^*\}$ be a semi-normalized $(1+\e_i)$-bounded M-basis.
  
   \Buo we may and do assume that $1\le \|x_i\|\le 1+\e_i$ and  
$\|x_i^*\|=1$ for all $i\in \bbN$ (if necessary, we may and do replace the given  M-basis by
 the equivalent $(1+\e_i)$-bounded M-basis defined by $y_i=\|x_i^*\|x_i$ and $y_i^*=x_i^*/\|x_i^*\|$).

Let $(e_i)_{i\in\bbN}$ be the standard orthonormal basis of $\ell_2$ and let $T:\ell_2\to\ell_2$ be the linear operator  defined by $Te_i=x_i$ for each  $i\in\bbN$. In order to prove that the operator $T$ is an isomorphism it is enough to prove that there exists a constant $K$ such that for each $n\in\bbN$ we have $ \|T_n\|\cdot\|T_n^{-1}\|\le K$, where $T_n$ is the restriction of the operator $T$ to 
 $\overline{\Span}\{e_i\}_{i=1}^n$.
 
Fix an arbitrary $n\in\bbN$ and consider the subset $\{x_i,x_i^*\}_{i=1}^n$ of the given M-basis. This subset is an M-basis of 
$\ell_2^n$. Thus,
 as described in Section~\ref{prelim}, it determines a unitary matrix $ U=[u_{ij}]\in {\mathcal{O}}(n)$
 and a positive diagonal matrix   $D=[d_i]$ 
  with $d_1\ge d_2\ge\dots\ge d_n>0$, that satisfy \eqref{bii}-\eqref{M-basis-matrix}. Since 
  $\e_i<1$ for all $i\in \bbN$, the matrices $U$ and $D$ satisfy 
\begin{align}\lb{epsilonj}
\forall 1\le i \le n\ \ \ \sum_{j=1}^n u_{ij}^2 d_j&\le (1+\e_i)^2< 1+3\e_i,\\
\lb{normalizedfunctionals}
\forall 1\le i \le n \ \ \ \sum_{j=1}^n u_{ij}^2 \frac{1}{d_j}&=1.
\end{align}
   
   To prove Theorem~\ref{(1)}, by \eqref{normofA}, it is enough  to prove that the ratio  $\sqrt{d_1}/\sqrt{d_n}$ is bounded by a constant independent of $n$.

   Let 
  $C\in \bbR$ be such  that $\sum_{i=1}^\infty \e_i \le C$.
  Then,  by \eqref{epsilonj} and \eqref{normalizedfunctionals}, we have
   \begin{align}
   \sum_{j=1}^n  d_j&\le n+3C,\\
   \sum_{j=1}^n  \frac{1}{d_j}&=n.
   \end{align}
 Hence  
\[
\sum_{j=1}^n \big(\sqrt{d_j}-\frac{1}{\sqrt{d_j}}\big)^2
=
\sum_{j=1}^n \big(d_j + \frac{1}{d_j}-2\big)
\leq (n+3C)+n-2n=3C,
\]
and therefore for each $j$ we have $\sqrt{d_j} \in [\frac{1}{R},R]$, where $R$ is the larger positive root of the equation 
$(x-\frac{1}{x})^2= 3C$. We thus obtain
$$
\frac{\sqrt{d_1}}{\sqrt{d_n}}\leq R^2 < R^2+\frac{1}{R^2} = 3C+2,
$$
which ends the  proof of Theorem~\ref{(1)}.

\section{Proof of Theorem~\ref{(2)}.}\lb{sec(2)}

Suppose that
 $(\e_i)_{i=1}^\infty$ is a nonnegative  sequence such that $\sum_{i=1}^\infty \e_i=\infty$. 
For $n\in \bbN$, we denote $s_n=\sum_{i=1}^n \e_i$ and $s_0=0$. We take a strictly increasing sequence 
$(n_m)_m\subset \bbN$, such that $s_{n_1}>1$, 
\begin{equation}\lb{infty}
\lim_{m\to \infty} (s_{n_m}-s_{n_{m-1}})=\infty,
\end{equation}
and the sequence $(r_m)_m$, where $r_m=\sqrt{s_{n_m}-s_{n_{m-1}}}$, is increasing.

We will construct a sequence of finite dimensional Hilbert spaces $(H_m)_{m\in\bbN}$ such that 
$\dim H_m={n_m}-{n_{m-1}}$, and an M-basis in $H_m$ whose  distance  to the orthonormal basis is equal to $r_m=\sqrt{s_{n_m}-s_{n_{m-1}}}$.

To construct $H_1$, take $n=n_1$  and  the vector $v=\frac{1}{\sqrt{s_n}}(\sqrt{\e_1},\dots,\sqrt{\e_n})$.
 Note that 
$\|v\|=1$. Let $U$ be an $n\times n$ orthogonal matrix whose first column vector is $ v$ and 
let $ D$ be a diagonal matrix $ D=[d_i]$, 
where 
 $d_1=s_n=r_1^2$ and 
  $d_i=1$ for $i=2,\dots,n$.

  Let  $B=UDU^T$.
  Then, for all $1\le i \le n$, we have
  \begin {align*}b_{ii}&=\sum_{j=1}^n u_{ij}^2 d_j
   = \frac{\e_i}{s_n}\cdot s_n +\sum_{j=2}^n u_{ij}^2=\e_i+(1-\frac{\e_i}{s_n})
  \in \Big(1, 1+\e_i\Big),\\
  (B^{-1})_{ii}&=\sum_{j=1}^n u_{ij}^2 \frac{1}{d_j}=\frac{\e_i}{s_n}\cdot \frac1{s_n} +\sum_{j=2}^n u_{ij}^2
 =\frac{\e_i}{s_n^2}+(1-\frac{\e_i}{s_n}) 
  \in \Big(1-\frac{\e_i}{s_n}, 1\Big).
  \end{align*}

 Thus, for all $1\le i \le n$,
\begin{equation}\lb{M-basis-2-matrix}
\ \ \  \ \ \ \left(\sum_j u_{ij}^2 d_j\right)\left(\sum_j u_{ij}^2\frac{1}{d_j}\right) < 1+\e_i< (1+\e_i)^2.
\end{equation}

Hence, by \eqref{M-basis-matrix},  the column vectors of the  the
  matrix $A=\sqrt{D}U^T$ form a  $(1+\e_i)$-bounded M-basis in $\ell_2^n$.
 (Note that, instead of taking the identity matrix, we could take  an arbitrary unitary matrix $V$ and define $A=V\sqrt{D}U^T$.)

 By \eqref{normofA}, the  distance 
 between this  M-basis and  the standard Riesz basis of $H_1=\ell_2^n$ is equal to 
 $\sqrt{d_1}/\sqrt{d_n}=\sqrt{s_n}=r_1$.

 Now we repeat the above argument replacing the sequence 
 $(\e_i)_{i=1}^\infty$ by $(\e_i)_{i=n+1}^\infty$ and we get the next finite dimensional Hilbert space $H_2$ with an $M$-basis whose distance to the orthonormal basis is equal to $r_2=\sqrt{s_{n_2}-s_n}$. Repeating this construction {\it ad infinitum} we get the sequence of finite dimensional Hilbert spaces $(H_m)_{m=1}^\infty$ each with an $M$-basis distant from the orthonormal basis by  $r_m=\sqrt{s_{n_m}-s_{n_{m-1}}}$, $m=1,2,\dots$, and by \eqref{infty}, the sequence $(r_m)_m$ diverges to infinity. Then, with a slight abuse of notation,  the union of all constructed $M$-bases of spaces $H_m$  forms the required $M$-basis in the space
 $\ell_2=\big(\bigoplus_{m=1}^\infty H_m\big)_{\ell_2}$.

\section{Proof of Theorem~\ref{(4)}.}
Let $\{\e_i\}_{i=1}^\infty$ be a nonnegative sequence satisfying \eqref{extra} and such that for all $i\in \bbN$, $\e_i< 1$. Since \eqref{extra} implies that $\sum_i\e_i=\infty$, we will use  a very similar construction of 
 the  M-basis $\{x_i\}_{i\in\bbN}$ as  in the proof of  Theorem~\ref{(2)}, but we will be a little more careful in our selection of the  sequence 
$(n_m)_m\subset \bbN$ (we  continue to use the same notations as in Section~\ref{sec(2)}), which besides satisfying  \eqref{infty} will  also  satisfy a stronger condition. To describe this condition, first observe  that it follows easily from the Stolz-Ces\`aro Theorem
(i.e. a discrete version of the L'Hopital's rule), that for any  nonnegative sequence $\{\e_i\}_{i=1}^\infty$  satisfying \eqref{extra}, that is, $\lim_n n\e_n=\infty$, for every $N\in \bbN$ we have  
 \begin{equation}\lb{extraold}
\lim_{n\to\infty}\frac{1}{\sqrt{n-N}}\sum_{i=N+1}^{n}\sqrt{\e_i}=\infty.
\end{equation}

Thus,
by \eqref{extraold},  there exists a  sequence 
$(n_m)_m\subset \bbN$ such that 
\begin{equation}\lb{new}
\forall m, \quad n_m-n_{m-1}\ge 32, \quad\text{and}\quad \lim_{m\to\infty}\frac{1}{\sqrt{n_m-n_{m-1}}}\sum_{i=n_{m-1}+1}^{n_m}\sqrt{\e_i}=\infty.
\end{equation}

We  construct an M-basis in the same way as in Section~\ref{sec(2)}, but using the sequence 
$(n_m)_m\subset \bbN$ satisfying \eqref{new}.
We will prove that no   permutation of  this M-basis  is  equivalent to a Schauder basis. 

\begin{remark} Note that we do not claim, that the permuted system 
  $\{x_{\s(i)}\}_{i\in\bbN}$ is a   $(1+\e_i)$-bounded M-basis in 
  $\ell_2$. Clearly  it is  a $(1+\e_{\s(i)})$-bounded M-basis in $\ell_2$, but  this is not important for us.
  \end{remark}

Let $\s$ be any permutation of $\bbN$. It is enough to show that for any  constant 
$C\ge 1$, there exist finite sets $E\subset F\subset \bbN$ satisfying $\nu<\la$ for every $\nu\in E$ and $\la\in F\setminus E$,  and coefficients $(\de_{\nu})_{\nu\in F}\subset \bbR$ such that 
\begin{equation}\lb{basisconst}
\frac{\Big\|\sum\limits_{\nu\in E} \de_{\nu} x_{\s(\nu)}\Big\|}{\Big\|\sum\limits_{\nu\in F} \de_{\nu} x_{\s(\nu)}\Big\|}\ge C.
\end{equation}

Let $C\ge 1$ be given. By \eqref{new}, we take $m\in \bbN$ such that 
\begin{equation}\lb{choicem}
\frac{1}{\sqrt{n_m-n_{m-1}}}\sum_{i=n_{m-1}+1}^{n_m}\sqrt{\e_i}\ge3C.
\end{equation}

Define 
  \[F=\s^{-1}(\{n_{m-1}+1,\dots,n_m\}),\]
  and let $F=\{k_j\}_{j={n_{m-1}+1}}^{n_m}$ be the enumeration of the elements of  the set $F$ in the strictly increasing order. 
Let 
$${ t_m=\sum_{i={n_{m-1}+1}}^{n_m} \sqrt{\e_i}=\sum_{\nu\in F} \sqrt{\e_{\s(\nu)}}=\sum_{j={n_{m-1}+1}}^{n_m} \sqrt{\e_{\s(k_j)}}}.$$

Note that the three expressions above define the same sum, however, in  the first one  the summation is in the   order of the  original M-basis $\{x_i\}_i$, while in the third  this finite sum is written in the order of the permuted basis  $\{x_{\s(i)}\}_i$\! .
Below we work with the permuted M-basis, so we use the order of the sum from the  third expression.

Let $\al$ be the unique number such that
\begin{equation}\lb{def-alm}
\sum_{j=n_{m-1}+1}^{\al} \sqrt{\e_{\s(k_j)}}\ge \frac{t_m}2\ \ \text{\rm \ and \ } \ 
\sum_{j=n_{m-1}+1}^{\al-1} \sqrt{\e_{\s(k_j)}}< \frac{t_m}2.
\end{equation}
Note that the number $\al$ depends both on $m$ and on the permutation $\s$.

We will prove that \eqref{basisconst} holds for the sets $F$ and 
$E=\{k_j\}_{j={n_{m-1}+1}}^{\al}\subset F,$
and the coefficients
\[\de_{\nu}=\begin{cases} 1 \ \ & \text{if }\ \nu\in E,\\
-1 \ \ & \text{if }\  \nu\in F\setminus E.
\end{cases}\]

Note that we have 
\[\Big\|\sum\limits_{\nu\in E} \de_{\nu} x_{\s(\nu)}\Big\|
=\Big\| \sum_{j=n_{m-1}+1}^{\al} \hde_j x_{\s(k_j)}\Big\|\ 
\text{\ and}\ \ 
\Big\|\sum\limits_{\nu\in F} \de_{\nu} x_{\s(\nu)}\Big\|=\Big\| \sum_{j=n_{m-1}+1}^{n_m} \hde_{j} x_{\s(k_j)}\Big\|,\]
where
\[\hde_{j}=\de_{k_j}=\begin{cases} 1 \ \ & \text{if }\ n_{m-1}+1\le j\le \al,\\
-1 \ \ & \text{if }\ \al+1\le j\le n_m.
\end{cases}\]

The set of vectors $\{x_{\s(\nu)}\}_{\nu\in F}$ is equal to the set $\{x_i\}_{i=n_{m-1}+1}^{n_m}$ which, by construction, is the $(1+\e_i)$-bounded  M-basis of the space $H_m$
and  consists  of the column vectors of the matrix 
$A_m=\sqrt{D_m}U_m^T$, where 
 $U_m=(u_{ij})_{i,j=n_{m-1}+1}^{n_m}$ is an  orthogonal matrix  with  the first column vector equal to $v^m=\frac{1}{\sqrt{s_{n_{m}}-s_{n_{m-1}}}}(\sqrt{\e_{i}})_{i=n_{m-1}+1}^{n_m}
 =\frac{1}{r_m}(\sqrt{\e_{i}})_{i=n_{m-1}+1}^{n_m}$ , 
 and 
 $ D_m$ is a diagonal matrix $ D=[d_i]_{i=n_{m-1}+1}^{n_m} $ with
  $d_{n_{m-1}+1}=s_{n_{m}}-s_{n_{m-1}}=r_m^2$ and 
  $d_i=1$ for $i=n_{m-1}+2,\dots,n_m$. 
  Thus, as in \eqref{bii}, if we denote  $B_m=U_m D_m U_m^T$, then for all $n_{m-1}<i,j\le n_m$, we have
  \begin{equation}\lb{bij}
\begin{split} 
  (x_{\s(k_i)}, x_{\s(k_j)})&=(B_m)_{\s(k_i)\s(k_j)}=\sum_{\nu=n_{m-1}+1}^{n_m} u_{{\s(k_i)}\nu}u_{{\s(k_j)}\nu}d_{\nu}\\
  &= u_{{\s(k_i)},{n_{m-1}+1}}u_{{\s(k_j)},{n_{m-1}+1}}d_{n_{m-1}+1}+
  \sum_{\nu={n_{m-1}+2}}^{n_m} u_{{\s(k_i)}\nu}u_{{\s(k_j)}\nu}d_{\nu}\\
   &=(r_m^2-1)  u_{{\s(k_i)},{n_{m-1}+1}}u_{{\s(k_j)},{n_{m-1}+1}}+
  \sum_{\nu={n_{m-1}+1}}^{n_m} u_{{\s(k_i)}\nu}u_{{\s(k_j)}\nu}\\
    &=\frac{(r_m^2-1)}{r_m^2}\sqrt{\e_{\s(k_i)}\e_{\s(k_j)}}+
     \sum_{\nu={n_{m-1}+1}}^{n_m} u_{{\s(k_i)}\nu}u_{{\s(k_j)}\nu}.
\end{split}
 \end{equation}
Since $U_m$ is orthogonal, the second summand in the last formula equals 1 if $i=j$, and 0 otherwise.

    By \eqref{bij},  for each $M$ with $n_{m-1}+1\le M\le n_m$ we have
   \begin{align}
   \Big\| \sum_{j=n_{m-1}+1}^M \hde_{j} x_{\s(k_j)}\Big\|^2
   &=\Big(\sum_{j=n_{m-1}+1}^M \hde_{j} x_{\s(k_j)}, \sum_{j=n_{m-1}+1}^M \hde_{j} x_{\s(k_j)}\Big)\notag\\
   &=\sum_{i,j=n_{m-1}+1}^M \hde_{i}\hde_j (x_{\s(k_i)}, x_{\s(k_j)})\notag\\
   &=\sum_{i,j=n_{m-1}+1}^M \hde_{i}\hde_j\Big[\frac{(r_m^2-1)}{r_m^2}
   \sqrt{\e_{\s(k_i)}\e_{\s(k_j)}}+
     \sum_{\nu={n_{m-1}+1}}^{n_m} u_{{\s(k_i)}\nu}u_{{\s(k_j)}\nu}\Big]\notag\\
   &=\sum_{i,j=n_{m-1}+1}^M \hde_{i}\hde_j\frac{(r_m^2-1)}{r_m^2}\sqrt{\e_{\s(k_i)}\e_{\s(k_j)}}
   +
     \sum_{i=n_{m-1}+1}^M \hde_{i}^2\notag\\
   \lb{norm} &=\frac{(r_m^2-1)}{r_m^2}\Big(\sum_{j=n_{m-1}+1}^M \hde_{j}\sqrt{\e_{\s(k_j)}}\Big)^2
   +
    \Big(M-n_{m-1}\Big).
   \end{align}

To estimate the first summand of \eqref{norm}, note that, by \eqref{def-alm}, we have
\[\frac{t_m}2 \le \sum_{j=n_{m-1}+1}^{\al} \sqrt{\e_{\s(k_j)}}< \frac{t_m}2+\sqrt{\e_{\s(k_\al)}}.\]
Therefore
\begin{equation}\label{small-estimate}
\begin{split}
\Big|\sum_{j=n_{m-1}+1}^{n_m} \hde_i\sqrt{\e_{\s(k_j)}}\Big|
&=\Big|\sum_{j=n_{m-1}+1}^{\al} \sqrt{\e_{\s(k_j)}}-\sum_{j=\al+1}^{n_m} \sqrt{\e_{\s(k_j)}}\Big|\\
&=\Big|2\Big(\sum_{j=n_{m-1}+1}^{\al}\sqrt{\e_{\s(k_j)}}\Big)-t_m \Big|< 2\sqrt{\e_{k_\al}}\le 2.
\end{split}
\end{equation}

  Thus, by \eqref{norm} with   $M=n_m$ 
  \remove{, \eqref{sum-of-unitary},} 
  and  \eqref{small-estimate}, we obtain
    \begin{equation}\lb{sum-all-B}
 \begin{split}
\Big\|\sum_{\nu\in F} \de_{\nu} x_{\s(\nu)}\Big\|^2&= \Big\| \sum_{j=n_{m-1}+1}^{n_m} \hde_{j} x_{\s(k_j)}\Big\|^2\\
 &=\frac{(r_m^2-1)}{r_m^2}\Big(\sum_{j=n_{m-1}+1}^{n_m} \hde_{j}\sqrt{\e_{\s(k_j)}}\Big)^2+n_m-n_{m-1}\\
&\le \frac{(r_m^2-1)}{r_m^2}\cdot2^2+n_m-n_{m-1}\le n_m-n_{m-1}+4.
   \end{split}
  \end{equation}

On the other hand,  by \eqref{norm} with  $M=\al$, \eqref{choicem}, and  
\eqref{def-alm}, 
  we obtain
    \begin{equation}\lb{sum-mn-B}
 \begin{split}
\Big\|\sum_{\nu\in E} \de_{\nu} x_{\s(\nu)}\Big\|^2&= \Big\| \sum_{j=n_{m-1}+1}^{\al} \hde_{j} x_{\s(k_j)}\Big\|^2\\
&=\frac{(r_m^2-1)}{r_m^2}\Big(\sum_{j=n_{m-1}+1}^{\al} \hde_{j}\sqrt{\e_{\s(k_j)}}\Big)^2
   + \Big(\al-n_{m-1}\Big)
   \\
&\ge  \frac12\Big(\frac{t_m}2\Big)^2  +0=\frac18 t_m^2.
   \end{split}
  \end{equation}
   
Combing \eqref{sum-all-B}, \eqref{sum-mn-B},  \eqref{choicem}, and, since by  \eqref{new},
 $ n_m-n_{m-1}\ge 32$,
 we obtain
\[ 
\frac{\Big\|\sum_{\nu\in E} \de_{\nu} x_{\s(\nu)}\Big\|}
{\Big\|\sum_{\nu\in F} \de_{\nu} x_{\s(\nu)}\Big\|}\ge 
 \frac{t_m}{\sqrt{8(n_m-n_{m-1}+4)}}\ge \frac{t_m}{\sqrt{9(n_m-n_{m-1})}}\ge C,
\]
which proves \eqref{basisconst}, and  thus ends the proof of Theorem~\ref{(4)}.

\section{Proof of Corollary~\ref{(5)}.}

We follow the same proof as in \cite{V99}, except that in the constuction in the proof of \cite[Lemma~2]{V99} in addition to the suspace $\widetilde{F}$ we take also another subspace $H$ such that $\dim H=\dim \widetilde{F}$ and 
$\widetilde{F}+H$ is 
$(1+\e)$ isomorphic to $\ell_2^{2(M-N)}$. We extend the M-basis to 
$\widetilde{F}$ exactly as it is done in \cite[Lemma~2]{V99}, and we take an 
M-basis of $H$ constructed as in Theorem~\ref{(4)}.

\section{Proof of Corollary~\ref{(3)}.}

This is the same proof as the proof of \cite[Corollary 1]{Pl77}. We include it here for the convenience of the reader.

Let $\e>0$. By Corollary~\ref{(5)}, there exists a $(1+\e)$-bounded M-basis  
$\{x_i,x_i^*\}_{i\in\bbN}$ in $X$ such that none of its permutations is    equivalent to  a Schauder basis of 
$X$ and such  that for all $i\in\bbN$,
  $\|x_i\|\cdot\|x_i^*\|\le 1+\e$,  $\|x_i\|=1$, and $\|x_i^*\|\le 1+\e$.  We define a new norm on 
  $X$ by
 \begin{equation*}\lb{normdef}
|\!|\!|x|\!|\!|\DEF\max\{\|x\|,\sup\{ |x_i^*(x)| : i\in\bbN\}\}.
\end{equation*} 

Then, for all  $i\in\bbN$, we have 
$ |\!|\!|x_i|\!|\!|=1$, and for all $x\in X$ 
\[ \|x\|\le |\!|\!|x|\!|\!| \le (1+\e)\|x\|.\]

Moreover, for all  $j\in\bbN$,
\[ |\!|\!|x_j^*|\!|\!|=\sup\{ |x_j^*(x)| : x \ \text{s.t.}\ |\!|\!|x|\!|\!|=1\}\le 1.\]
Since $x_j^*(x_j)=1$, we conclude that $|\!|\!|x_j^*|\!|\!|=1$, and thus
 $\{x_i,x_i^*\}_{i\in\bbN}$ is an Auerbach basis in 
$(X, |\!|\!|\cdot|\!|\!|)$.

\section{Examples}\lb{examples}

\begin{example}\lb{exl1}
We begin with showing that in $\ell_1$ there exists an Auerbach basis that is not equivalent to a Schauder basis under any permutation.
\end{example}

 \begin{proof}[Outline of proof]
 Let $(G_n)_{n=1}^\infty$ be a sequence of finite abelian groups
with ranks tending to infinity and equal to powers of 2.
Clearly $\big(\bigoplus_n L^1(G_n)\big)_1$ is isometrically isomorphic to $\ell_1$ (here $L^1$ is taken with respect to the normalized Haar measures). It is also clear that the system 
\[
\{(0,\dots,0,\vertarrowbox{\chi}{\text {\small\  $n$-th}},0,\dots) : \text{$\chi$ is a character of $G_n$}\}
\] 
is an Auerbach basis (dual functionals are the same characters just treated as elements of $L^\infty(G_n)$).
Suppose that in some permutation this system is a Schauder basis with
basis constant $C>0$. Let $n\in \bbN$ be sufficiently big and $\chi_1,\chi_2,\dots,\chi_{2^m}$ be the enumeration of characters of $G_n$ in order compatible with this permutation. Clearly, for every $j=1,2,\dots,2^m$, we get that $\nu_j=\chi_1+\chi_2+\dots+\chi_j$ is an idempotent measure with norm not exceeding $C$. The support of the Fourier transform of $\nu_j$, denoted by $A_j$, has exactly $j$ elements. By the quantitative version of Helson's Idempotent Theorem  due to Green and Sanders (see \cite{gs,Sa}), for every $j=1,2,\dots, 2^m$, we have
$$
{\bf 1}_{A_j}=\sum_{i=1}^L \varepsilon_i {\bf 1}_{K_i},
$$
where $K_i$ are cosets of $G_n$, $\varepsilon_i$ are signs, and $L< \exp\left( C^{4+o(1)}\right)$ (for us it is only important that $L$ is a fixed number, independent of $n$). Counting the number of elements we get immediately
$$
\# A_j=\sum_{i=1}^L \varepsilon_j \cdot \# K_i
$$
Observe now that the number of elements of any coset of $G_n$ is a power of 2. Thus each of the numbers 
$\varepsilon_i\cdot \#K_i$ can take only $2m$ possible values, so the whole sum can take only $(2m)^L$ possible values.
But for big $m$ obviously $2^m>(2m)^L$ and we get a contradiction.
\end{proof}

\begin{example}\lb{exlp}
For any  $1<p<\infty$, $p\neq 2$, the space $\ell_p$ contains an Auerbach basis that is 
not an unconditional Schauder basis.
\end{example}
 
 \begin{proof}[Outline of proof]
\remove{Next example is the space $\ell_p$ for $p\neq 2$ for which we show the existence of an Auerbach basis which is not unconditional Schauder basis.}
Clearly $\ell_p=\big(\bigoplus_n L^p(G_n)\big)_p$. This time we just want the rank of $G_n$ to tend to infinity, without any special arithmetical requirements. The Auerbach system is the same as in the first example. If this system was an unconditional basis, we would have
$$
\aligned
\|\sum a_j\chi_j\|^p_{L^p(G_n)} &\simeq
\int_0^1\|\sum r_j(t)a_j\chi_j\|^p_{L^p(G_n)}
=\int_0^1\int_{G_n}|\sum r_j(t)a_j\chi_j|^p\,dtds\\
&\simeq \int_0^1\int_{G_n}(\sum |a_j|^2)^{p/2}\,dtds
=(\sum |a_j|^2)^{p/2}
\endaligned $$
Hence by the Parseval's theorem $L^p(G_n)$ would be uniformly isomorphic to $L^2(G_n)$.

Note that a very similar argument works also for the spaces
$\big( \bigoplus_n L^{p_n}(G_n)\big)_2$ constructed by Johnson as examples of HAPPY spaces (cf. \cite{J79}).
  \end{proof}

Example~\ref{exlp}  can be extended to wider class of 1-symmetric sequence spaces.

\begin{example} \lb{exsym}
 Let $X$ be a Banach space of type 2 or  of cotype 2 (see \cite[Definition~1.e.12]{LT2}) with
a 1-symmetric basis $(e_n)_{n=1}^\infty$ such that $X$ is not isomorphic to a Hilbert space. Then there exists an Auerbach basis in $X$ that is not an unconditional Schauder basis.

 In particular, this holds for Orlicz sequence  spaces not isomorphic to a Hilbert space, satisfying the $\Delta_2$ condition at zero, with the Matuszewska-Orlicz lower index greater than or equal to 2 (or, by duality,  the Matuszewska-Orlicz upper index smaller than or equal to 2), cf. \cite{KMP96}.
\end{example}

 \begin{proof}[Outline of proof]
Suppose $X$ has type 2.  If $X$ has cotype 2, then we pass to its dual.
Let $(n_k)_{k=1}^\infty$ be an increasing sequence of natural numbers and $X_k= \Span \{e_j:  n_1+\dots+n_{k-1}<j\leq n_1+\dots n_{k-1}+n_k\}$. Clearly, $X_k$ is  isometrically lattice isomorphic to some space $E(\bbZ_{\dim X_k})$ of functions on the group $\bbZ_{\dim X_k}$ with the norm invariant with respect to the group action. Let us denote this isomorphism by
$T_k:E(\bbZ_{\dim X_k})\to X_k$. Let $F_k=\{\chi^k_j: 1\leq j\leq \dim X_k\}$ be the set of all characters of $\bbZ_{\dim X_k}$. Then it is easy to see that the set
$F=\bigcup_{k=1}^\infty T_k(F_k)$ is an Auerbach basis in $X$
(the biorthogonal functional corresponding to $T_k(\chi^k_j)$ is $(T_k^\ast)^{-1}(\chi^k_j)$). 

If the set $F$ is an unconditional Schauder basis, then an
 argument similar to the one in  Example~\ref{exlp}  shows that spaces
$X_k=\Span \{T_k(\chi^k_j): 1\leq j\leq \dim X_k\}$ are uniformly isomorphic to Hilbert spaces. Indeed, by the Jensen inequality, 
$$
 \int_0^1 \|\sum r_j(t)a_j\chi^k_j\| \,dt \geq 
 \|\int_0^1|\sum r_j(t)a_j\chi^k_j|\,dt \|.\remove{\simeq(\sum |a_j|^2)^{1/2}}
$$
And, by Khinchin inequality and since the characters attain only values of modulus 1, at any point $s\in \bbZ_{\dim X_k}$ we have
$$
\int_0^1 |\sum r_j(t)a_j\chi^k_j|\simeq (\sum |a_j|^2)^{1/2}.
$$
On the other hand, since $X$ has type 2, and thus
$$
\int_0^1 \|\sum r_j(t)a_j\chi^k_j\| \,dt \leq (\sum \|a_j\chi^k_j\|^2)^{1/2}\simeq(\sum |a_j|^2)^{1/2}.
$$
By unconditionality, $\|\sum a_j\chi^k_j\|$ is uniformly proportional to $(\sum |a_j|^2)^{1/2}$, which implies 
 that spaces $X_k$ are uniformly isomorphic to Hilbert spaces.
 
 Since $X$ is symmetric and $\dim(X_k)=n_k$ increases to $\infty$, we conclude  that   $X\simeq \ell_2$, contrary to our assumption.
\end{proof}

The next example combines Examples~\ref{exl1} and \ref{exsym}. 

\begin{example}\lb{ex4}
If   a space $X$ with a 1-symmetric basis $(e_n)_{n=1}^\infty$   is
``sufficiently close'' to $\ell_1$, then there exists an Auerbach basis in $X$ that is not a Schauder basis under any permutation.
\end{example}

 \begin{proof}[Outline of proof]
The construction is exactly the same as in Example~\ref{exsym} with $n_k=2^k$ for $k=1,2,\dots$. Let $S_k: X_k\to\ell_1^{2^k}= L^1(\bbZ_{2^k})\subset\ell_1$ be the linear extension of the operator given by $S_k(e_j)= h_j$, for $j$ with
 $n_1+n_2+\dots+n_{k_1}<j\leq n_1+n_2+\dots +n_k$, where $h_j$ are elements of the standard basis of $\ell_1$. Let $C_k=\|S_k\|$ (note that  $\|S_k^{-1}\|=1$ since we assume that the 1-symmetric basis 
 $(e_j)_j$ of $X$ is normalized).
  We measure the closeness of $X$ to $\ell_1$ by the rate of growth of
 $C_k$.
 
  The operator $S_k$ transforms $X_k=E(\bbZ_{2^k})$ onto the set of characters of $\bbZ_{2^k}$, and to any bounded  projection $P$ from $X_k$ onto a subspace of $X_k$ spanned by a subset of characters corresponds the projection $S_kPS_k^{-1}$ of $L^1(\bbZ_{2^k})$
 onto the subspace spanned by the corresponding set of characters.
 As in Example~\ref{exl1} we derive that if any permutation of characters is a Schauder basis with  basis constant $C$, then there has to be $2^k$ of such projections of pairwise different ranks, and the Green-Sanders theorem forces that
 $2^k< (2k)^L$ where 
 $L= \max_P\big(\exp (\|S_kPS_k^{-1}\|^{4+o(1)})\big)\le \exp (C\cdot C_k^{4+o(1)})$. This leads to a contradiction if $C_k<(\log k)^{1/4-\varepsilon}$ for any $\varepsilon>0$ (remember that 
 $C_k$ is the value of the norm for a $2^k$-dimensional space spanned by the elements of the basis; this is a double logarithm of the actual dimension).
\end{proof}

\begin{example}\lb{expel}
The last example is conditional. Let us recall a conjecture of   Pe\l czy\'nski that says that if $p\neq 2$, then the characters  in the space
$L^p({\bf T}^\infty)$ do not form a Schauder basis in any permutation 
\cite[Problem 5.1]{Pe2006}. Since, obviously, the characters do form an Auerbach basis, then the validity of the  Pe\l czy\'nski conjecture would give an example of an Auerbach basis which is not a Schauder basis under any permutation  in the space $L^p$. 
\end{example}


\begin{small}

\renewcommand{\refname}{
\end{small}

\textsc{Department of Mathematics, Miami University, Oxford, OH
45056, USA} 

  \textit{E-mail address}: \texttt{randrib@miamioh.edu} \par

\textsc{Institute of Mathematics, Polish Academy of Sciences, ul. Śniadeckich 8, 00-656 Warsaw, Poland}

  \textit{E-mail address}: \texttt{miwoj@impan.pl} \par
  
  \textsc{Department of Mathematical Sciences, University of Cincinnati, Cincinnati, OH
45221, USA}

\textit{E-mail address}: \texttt{zatitspl@ucmail.uc.edu}

\end{document}